\documentclass{amsart}

\usepackage{amsmath,amssymb}
 
\theoremstyle{thmit} 
\newtheorem{thm}{Theorem}[section]
\newtheorem{lem}[thm]{Lemma}
\newtheorem{cor}[thm]{Corollary}

\newtheorem{obs}[thm]{Observation}
\newtheorem{EXA}[thm]{Example}

\theoremstyle{thmrm} 
\newtheorem{exa}{Example}[section]

\newtheorem*{rem}{Remark}
\newtheorem*{oldproof}{Proof}

 \renewcommand{\leq}{\leqslant}
 \renewcommand{\geq}{\geqslant}
 \newcommand{\notdiv}{\mid{\kern -7.4pt/}\kern2pt}
 \renewcommand{\pmod}[1]{\mkern6mu({\rm mod}\mkern3mu#1)}
 \newcommand{\Ints}{\mathbb Z}
 \newcommand{\Nats}{\mathbb N}
 \newcommand{\Rats}{\mathbb Q}
 
 \newcommand{\Complex}{\mathbb C}
 \newcommand{\Field}{\mathbb F}
 \newcommand{\calK}{\mathcal{K}}
 \newcommand{\calC}{\mathcal{C}}

 \newcommand{\Sym}{\mathrm{Sym}\kern1pt}
 \newcommand{\Alt}{\mathrm{Alt}\kern1pt}
 \newcommand{\GL}{\mathrm{GL}\kern1pt}
 \newcommand{\AGL}{\mathrm{AGL}\kern1pt}
 \newcommand{\ASL}{\mathrm{ASL}\kern1pt}
 \newcommand{\SL}{\mathrm{SL}\kern1pt}
 \newcommand{\PSL}{\mathrm{PSL}\kern1pt}
 \newcommand{\SO}{\mathrm{SO}\kern1pt}
 \newcommand{\End}{\mathrm{End}\kern1pt}
 \newcommand{\Hom}{\mathrm{Hom}\kern1pt}
 \newcommand{\Char}{\mathrm{char}\kern1pt}
 \newcommand{\Aut}{\mathrm{Aut}\kern1pt}
 \newcommand{\Soc}{\mathrm{Soc}\kern1pt}
 \newcommand{\minN}{\hbox{\rm min-{\sc n}}}
 \newcommand{\minSN}{\hbox{\rm min-{\sc sn}}}
 \newcommand{\wreath}{\,\mathrm{wr}\,}
 \newcommand{\K}{\kern1pt}
 \newcommand{\Mid}{\kern3.5pt
              \rule[-3pt]{0.35pt}{12.6pt}\kern4pt}
 \newcommand{\Normal}{\kern2pt{\leq}\kern-1.3pt
              \rule[0.85pt]{0.43pt}{6.2pt}\kern3.5pt}

 \title{Some infinite permutation groups and 
 related finite linear groups}
 \author{Peter M.\ Neumann}
 \address{Peter M.\ Neumann\\The Queen's College\\
 Oxford OX1 4AW\\
 United Kingdom}
 \author{Cheryl E.\ Praeger}
 \address{Cheryl E.\ Praeger\\School of Mathematics and Statistics\\
 The University of Western Australia\\
 35 Stirling Highway, Crawley WA 6009\\ 
 Australia}
 \author{Simon M.\ Smith}
 \address{Simon M.\ Smith\\Department of Mathematics\\
 NYC College of Technology\\
 City University of New York (CUNY)\\
 NY, USA}


\begin{document}

 \maketitle

 \begin{center}
 \textit{Dedicated to the memory of our late colleague 
 and friend \\ L.\,G.\ {\rm(}Laci\/{\rm)} Kov\'acs}
 \end{center}

 \begin{abstract}
 This article began as a study of the structure of 
 infinite permutation groups $G$ in which point stabilisers 
 are finite and all infinite normal subgroups are transitive.
 That led to two variations.
 One is the generalisation in which point stabilisers 
 are merely assumed to satisfy \minN, the minimal condition 
 on normal subgroups. 
 The groups $G$ are then of two kinds.
 Either they have a maximal finite normal subgroup, 
 modulo which they have either one or two minimal 
 non-trivial normal subgroups, or they have
 a regular normal subgroup $M$ 
 which is a divisible abelian $p$-group of finite rank. 
 In the latter case the point stabilisers are finite 
 and act irreducibly on a $p$-adic vector space 
 associated with~$M$.
 This leads to our second variation, which is a study 
 of the finite linear groups that can arise.
 \end{abstract}

 \section{Introduction}\label{s:1}

 Stimulated by the O'Nan--Scott theory described in~\cite{SS2} 
 of primitive permutation groups that have 
 finite point stabilisers, we initiated a study of 
 infinite permutation groups in which stabilisers are finite 
 and all infinite normal subgroups are transitive.
 This class includes all primitive, or more generally 
 quasiprimitive, groups with finite point-stabilisers.
 Although infinite permutation groups with finite stabilisers 
 arise naturally in various 
 contexts they do not usually have the property that 
 their infinite normal subgroups are transitive.
 A crystallographic group, for example, has finite 
 stabilisers (point groups), but most of its infinite 
 normal subgroups are not transitive on its point-orbits.
 However, if an infinite permutation group~$G$ is primitive 
 (or even if it is no more than quasiprimitive), 
 then a point stabiliser $G_\alpha$ is finite if and only if 
 there is a finite upper bound on the lengths of the 
 $G_\alpha$-orbits%
 ---this is a special case of a theorem proved by 
 Schlichting \cite{Sch} 
 and independently by Bergman and Lenstra \cite{BL} that gives 
 necessary and sufficient conditions for a transitive group 
 to have a bound on its subdegrees, that is on the 
 lengths of orbits of a point stabiliser. 

 It was something of a surprise to us that our ideas about 
 groups in which all infinite normal subgroups are transitive 
 and stabilisers are finite could be naturally generalised to 
 those in which the stabilisers merely satisfy \minN, 
 the minimal condition on normal subgroups.
 (Philip Hall introduced the notation min-$n$, but $n$
 has too many other natural meanings in our mathematics, 
 so we use a variant.)

 \bigskip
 {\sc Notation, assumptions and terminology}.
 \begin{itemize}
 \item
 Throughout this paper $\Omega$ denotes an infinite set, $G$ 
 denotes a subgroup of $\Sym(\Omega)$ with the property that 
 all its infinite normal subgroups are transitive, 
 and $H := G_\alpha$, the stabiliser of $\alpha$, where 
 $\alpha \in \Omega$.
 \item
 We assume that $H$ satisfies \minN.
 \item
 If all non-trivial normal subgroups of a group $X$ are
 infinite (equivalently, if $\{1\}$ is the maximal finite 
 normal subgroup of $X$) then, for want of a better term, 
 we shall say that $X$ is {\it normally infinite}. 
 \end{itemize}

 Note that any quasiprimitive group of permutations of 
 an infinite set is normally infinite since non-trivial 
 normal subgroups, being transitive, are infinite. 

 To provide context, here are some 
 simply described, but in some sense representative, 
 examples of groups $G$ satisfying our conditions.

 \begin{exa}\label{Ex:affine} 
 Let $F$ be an infinite field, let $H := \SL(2,F)$, 
 and let $V := F^2$ with the natural action of $H$.
 Take $\Omega := V$ and $G := \ASL(2, F)$, the split 
 extension of the translation group of $V$ by $H$.
 Here $H$ is the stabiliser of $0$ and satisfies 
 \minN\ (it has centre of order $\leq 2$, 
 modulo which it is simple).
 The translation group is the unique minimal normal 
 subgroup.
 In this case $G$ is doubly transitive.
 \end{exa}

 \begin{exa}\label{Ex:Monolithic}
 Let $G$ be a simple group acting transitively
 on an infinite set $\Omega$ such that a stabiliser has \minN,
 for example, a stabiliser is finite. 
 Or let  $G := T \wreath_\Gamma H$ where $T$ is an infinite 
 simple group and $H$ is finite acting faithfully and 
 transitively on a set $\Gamma$ and $\Omega := T^\Gamma$.
 Here $T^\Gamma$, the base group of the wreath product, is  
 the unique minimal normal subgroup of $G$ and acts regularly, 
 and $H$ is a stabiliser.
 \end{exa}
 
 \begin{exa}\label{Ex:Bilithic}
 For any infinite simple group $T$, let $\Omega := T$,
 and let $G := T \times T$ acting by left and right 
 multiplication on $\Omega$ 
 (that is, $\omega^{(a,b)} = a^{-1}\omega b$).
 This  has two regular minimal normal subgroups, 
 each isomorphic to $T$, and the stabiliser $H$ of $1$ 
 is the diagonal.
 Then $H \cong T$, so obviously $H$ satisfies \minN.
 \end{exa}

 \begin{exa}\label{Ex:Pruefer} 
 For a prime number $p$ let $C_{p^\infty}$ 
 denote the Pr\"ufer \hbox{$p$-group} 
 (isomorphic to $\{\theta \in \Complex \Mid \exists k \in \Nats: 
 \theta^{p^k} = 1\} \leq \Complex^\times$). 
 If $\Omega := G := C_{p^\infty}$ with the regular action,
 then $G$ has only one infinite normal subgroup, namely $G$ itself,
 but arbitrarily large finite normal subgroups.   
 \end{exa}
 
 It will be convenient to have some terminology
 for phenomena illustrated in very basic form by 
 these examples.

 \begin{itemize}
  \item 
  A normally infinite permutation group that has 
  an abelian regular minimal normal subgroup 
  (as in Example~\ref{Ex:affine}) will be said 
  to be of {\it affine type}.
  \item 
  A normally infinite permutation group that has 
  a unique minimal normal subgroup that is non-abelian 
  (as in Example~\ref{Ex:Monolithic}) 
  will be said to be of {\it monolithic type}.
  \item 
  A normally infinite permutation group that has 
  precisely two minimal normal subgroups (each
  of which necessarily acts regularly, as in
  Example~\ref{Ex:Bilithic}) 
  will be said to be of {\it bilithic type}.
  \item 
  If, for some prime number $p$,
  our group $G$ has a regular normal subgroup that is
  a divisible abelian $p$-group of finite rank (hence is a 
  direct sum of finitely many copies of $C_{p^\infty}$---%
  as in Example~\ref{Ex:Pruefer}) 
  then $G$ will be said to be of 
  {\it $p$-divisible affine type}.   
 \end{itemize}

 Before stating our main theorems (to be
 proved in later sections), we give a further item 
 of contextual information. 

 \begin{obs}\label{Obs:GminN}
 Under our assumptions, $G$ satisfies \minN.
 \end{obs}

 {\it Proof}.
 Let $\mathcal N$ be any non-empty set of normal subgroups of~$G$.
 We show that $\mathcal N$ has minimal members.
 If $\mathcal N$ contains any finite normal subgroups 
 of $G$ then it contains one of smallest order, 
 and clearly this is minimal.
 Suppose now, therefore, that all members of $\mathcal N$ 
 are infinite.
 By the assumption on $G$, they are transitive on $\Omega$.
 Define 
 $\mathcal{N}_\alpha := \{ N \cap H \mid N \in \mathcal{N}\,\}$.
 Since $H$ satisfies {\minN} and all members of $\mathcal{N}_\alpha$ 
 are normal subgroups of $H$, there exists 
 $N_0 \in \mathcal N$ such that $N_0 \cap H$ is 
 minimal in~$\mathcal{N}_\alpha\,$.
 Suppose that $N \in \mathcal{N}$ and $N \leq N_0$.
 Then $N \cap H  = N_0 \cap H$ since $N_0 \cap H$ is minimal 
 in $\mathcal{N}_\alpha$ and $N \cap H \leq N_0 \cap H$.
 Now if $x \in N_0$ then since $N$ is transitive 
 on $\Omega$ there exists $y \in N$ such that 
 $\alpha\K y = \alpha \K x$, and so 
 $x = (xy^{-1})y \in (N_0 \cap H).N$, whence 
 (since $N_0 \cap H = N \cap H$), $x \in N$.
 Thus $N_0 = N$ and we have shown that $N_0$ 
 is minimal in $\mathcal{N}$.
 Hence $G$ satisfies \minN. 

 \bigskip
 Note that the Axiom of Choice (AC) is not needed in 
 the above proof.
 In fact, there are, we believe, only a few places where 
 it is really needed (in some cases in a weak form) 
 in this paper.
 Those will be noted.
 
 \bigskip
 Clearly, in any group $X$, either there are arbitrarily 
 large finite normal subgroups or there is a bound on 
 the sizes of finite normal  subgroups. 
 In the latter case, since the product of two finite 
 normal subgroups is a finite normal subgroup 
 there will be a unique maximal (largest) finite normal 
 subgroup $K$ and $X/K$ is normally infinite.
 
 \begin{thm}\label{thm:maxfin} 
 Suppose that $G$ as specified above has a maximal finite 
 normal subgroup $K$. 
 Then $K$ is semi-regular on $\Omega$ {\rm(}stabilisers 
 $K_\omega$ are trivial for all $\omega \in \Omega${\rm)}.
 If $\bar{G} := G/K$,  $\bar H := HK/K \cong H$, 
 and $\bar{\Omega} := \Omega/K = \Omega/\rho$ where $\rho$ 
 is the $G$-congruence whose blocks are the $K$-orbits,
 then $\bar G$ acts faithfully as a normally infinite group 
 on $\bar \Omega$ with stabiliser $\bar H$. 

 \medskip
 Moreover, if $G$ is normally infinite 
 {\rm(}equivalently, if $G$ is quasiprimitive 
 on $\Omega${\rm)} then {\bf either}
 $G$ has a unique minimal normal subgroup $M$  
 {\bf or} $G$ has {\rm(}precisely{\rm)} two minimal normal 
 subgroups $M_1,\, M_2$.

 If $M$ is unique either it is
 abelian and regular {\rm(}so $G$ is of affine type{\rm)} 
 or it is nonabelian {\rm(}so $G$ is of monolithic type{\rm)}
 and \hbox{$C_{G}(M) = \{1\}$}. 

 In the case of two minimal normal subgroups 
 {\rm(}where $G$ is of bilithic type{\rm)} 
 each $M_i$ acts regularly on $\Omega$, and 
 the group that they generate is their direct product.
 Also, if $M_0 := H\cap (M_1\times M_2)$ then $M_0$ 
 is a minimal normal subgroup of $H$ and projects isomorphically
 onto each of $M_1$ and $M_2$ {\rm(}hence $M_1 \cong M_2$ and
 $M_0$ is the diagonal of the direct product{\rm)}.
 \end{thm}

 This description of the possibilities in the case that $G$
 has a maximal finite normal subgroup probably cannot be 
 developed much further in general.
 When $G$ is normally infinite and monolithic $H$ acts 
 faithfully by conjugation as a group of automorphisms 
 of $M$, and all we know about $M$ is that it is  
 characteristically simple.
 One possibility is that it is a direct product of isomorphic 
 simple groups.
 Then its simple direct factors are minimal normal subgroups 
 of $M$ and they are permuted transitively under the 
 conjugation action of $H$.
 In this case $G$ is a sort of wreath product (perhaps twisted)
 of a simple group $T$ by $H$.
 Other possibilities are that $M$ could be a variant of 
 the McLain group (see \cite{DG, McL}) or one of Philip Hall's
 wreath powers \cite{Hall}, where in each case the relevant 
 index set is a dense linear ordering whose automorphism group
 contains a subgroup isomorphic with $H$ having an orbit
 that is unbounded both above and below. 
 For example, the index set could be 
 $\Rats$ with $H = \Aut(\Rats, \leq)$,
 a group that certainly satisfies \minN.
 These are just a few possibilities---it seems probable 
 that there are very many more.
  
 Very similar remarks apply to the bilithic case.
 Since $M_0$ is a minimal normal subgroup of $H$ it
 is characteristically simple, and any 
 characteristically simple group that can serve as 
 the socle of a monolithic group $G$ could serve as
 one of the two minimal normal subgroups of a group $G$
 of bilithic type.
 If we strengthen the condition on $H$ and suppose that it 
 satisfies \minSN, the minimal condition on subnormal 
 subgroups (clearly much stronger than \minN),
 then $M_0$ will have minimal normal subgroups.
 It then follows that $M_0$ is a direct product of isomorphic
 simple groups $T_i$ permuted transitively under conjugation
 by $H$, hence under~$H/M$.
 Thus in this case the minimal normal subgroups $M_1,\, M_2$
 of $G$ will also be direct products of simple groups,
 the simple factors in each being permuted transitively
 under conjugation by~$H$ 
 (see Observation~\ref{Obs:minSN} below).
 In particular, a little more can be said when $H$ is finite.

 \begin{thm}\label{thm:finitestabs}
 With the notation and assumptions specified above, 
 if $H$ is finite and $G$ is normally infinite then 
 $G$ is of monolithic type and its monolith $M$ 
 is a direct product $T_1 \times \cdots \times T_q$ of 
 finitely many isomorphic infinite simple groups.
 \end{thm} 
     
 If $G$ does not have a maximal finite normal subgroup then
 its structure is very different.

 \begin{thm}\label{ThmAffine}
 Suppose that $G$ {\rm(}as specified above{\rm)} has arbitrarily 
 large finite normal subgroups.
 Then $G$ has a unique minimal infinite normal\break 
 subgroup $M$,
 which acts regularly on $\Omega$. 
 For some prime number $p$,\break 
 $M$ is a divisible abelian $p$-group of finite rank\/ 
 {\rm(}so $G$ is of $p$-divisible affine type\/{\rm)}.
 Moreover, $H$ is finite and acts faithfully 
 and $p$-adic irreducibly 
 {\rm(}in the sense of 
 Theorem~{\rm\ref{Thm:H p-adic irreducible}} 
 below\/{\rm)} by conjugation on $M$.
 \end{thm}

 {\sc Remark~1}.\quad
 The proof that the rank of $M$ is finite requires AC.
 When $G$ is of this type, $\Omega$ and $G$ are 
 countably infinite (this also requires AC).
 Thus if $G$ is uncountable then it {\it must} have a 
 maximal finite normal subgroup.
  
 \bigskip
 {\sc Remark~2}.\quad
 If $H$ is finite then $G$ is either of $p$-divisible 
 affine type or it is almost monolithic,
 that is, an extension of a finite normal subgroup 
 acting semi-regularly by a twisted wreath product of 
 an infinite simple group by a finite group.
 This is an immediate consequence of 
 Theorems~\ref{thm:finitestabs} and~\ref{ThmAffine}.
 
 \medskip
 Proofs of Theorems~\ref{thm:maxfin} and~\ref{thm:finitestabs} 
 are given in~\S\ref{s:2} and a proof of Theorem~\ref{ThmAffine}
 is given in~\S\ref{s:3}.

 \medskip
 Theorem~\ref{ThmAffine} leads to an interesting question
 about finite groups and their modular representation theory.
 Suppose that $G$ is of $p$-divisible affine type, 
 so that $H$ is finite and $M$ is a divisible 
 abelian $p$-group of finite rank $r$,
 where $p$ is a prime number. 
 Let $V := M[p]$, the elementary abelian $p$-group
 $\{x \in M \mid x^p = 1\}$, 
 so that $V$ may naturally be construed
 as an $\Field_pH$-module of dimension $r$.
 A question that naturally arises is: 
 what pairs $(H, V)$ can occur?  
 In other words, what finite (linear) 
 groups can act faithfully (and $p$-adic
 irreducibly) on divisible abelian $p$-groups 
 of finite rank?
 We shall show that if $p$ is odd then $H$ must act
 faithfully on $V$, while if $p = 2$ then the kernel of 
 the action of $H$ is an elementary abelian $2$-group
 (Theorem~\ref{Th:faithful}).
 Moreover, a given linear group $(H, V)$ can arise 
 from a faithful action on a 
 divisible $p$-group if and only if $V$ is `liftable' 
 to an integral representation of~$H$ 
 (Theorem~\ref{Th:integralcover}).

 \section{Proofs of Theorems~\ref{thm:maxfin} 
 and~\ref{thm:finitestabs}}\label{s:2}

 Recall the notation and conventions from near the
 beginning of \S\ref{s:1}.
 In this section we assume that there is a bound on the 
 sizes of finite normal subgroups of $G$, so that there is a
 maximal finite normal subgroup $K$ and $G/K$ is normally
 infinite.

 Since $K$ is finite, the group of automorphisms of $K$ 
 induced by the conjugation action of $G$ is finite
 and therefore $C_G(K)$ is a normal subgroup of finite index 
 in $G$, hence infinite and so transitive on $\Omega$. 
 Then since $K$ has a transitive centraliser, 
 $K$ acts semi-regularly.
 (This is standard---here is the reason.
 If $x \in K_\omega$ and $\omega' \in \Omega$ 
 then there exists $y \in C_G(K)$ such that 
 $\omega' = \omega^y$; then 
 $(\omega')^x = (\omega')^{y^{-1}xy} = \omega^{xy} 
 = \omega^y = \omega'$, so $x$ fixes every point 
 of $\Omega$, whence $x = 1$.)

 In the statement of Theorem~\ref{thm:maxfin} we defined
 $\rho$ to be the equivalence relation on $\Omega$ whose
 classes are the $K$-orbits.
 Since $K \Normal G$, $\rho$ is a $G$-congruence 
 and $G$ acts transitively on $\Omega/\rho$. 
 If $L$ is the kernel of this action then $K \leq L$. 
 If $L$ were not equal to $K$ then, by the maximality of $K$, 
 $L$ would be infinite and hence transitive on $\Omega$, 
 which is not the case since $L$ has the same orbits as $K$. 
 Hence the kernel of the $G$-action on $\Omega/\rho$ is $K$.
 We define $\bar{G} := G/K$, $\bar{\Omega} := \Omega/\rho$,
 $\bar{\alpha} := \alpha K$ and $\bar{H} := HK/K$.
 With this notation, $\bar G$ acts faithfully
 on $\bar \Omega$, $\bar{G}_{\bar{\alpha}} = \bar{H}$,
 all infinite normal subgroups of $\bar G$ are transitive, 
 since $\bar{H} \cong H$ the stabiliser $\bar H$ satisfies 
 \minN, and since $K$ is the maximum finite normal 
 subgroup of~$G$, $\bar G$ is normally infinite.

 Now assume that $K = \{1\}$, so that $G$ is normally infinite.
 Since, as was shown in Observation~\ref{Obs:GminN}, 
 $G$ satisfies \minN, there are minimal (non-trivial) 
 normal subgroups in $G$. 
 Suppose first that there is just one minimal normal subgroup~$M$.
 Then $M$ is transitive so $G=MH$. 
 Since $C_H(M)$ is normalised both by $H$ and by $M$ it is normal
 in $G$, so it is trivial since $H$ contains no non-trivial
 normal subgroup of $G$.
 Therefore $C_H(M) = \{1\}$, so $C_G(M) = M$ or $C_G(M) = \{1\}$. 
 If  $C_G(M) = M$  then $M$ is abelian and 
 $H \cap M \leq C_H(M) = \{ 1 \}$,
 and so $M$ acts regularly and $G$ is of affine type.
 If $C_G(M) = \{1\}$ then $M$ is non-abelian and $G$ 
 is of monolithic type. 
 
 Now suppose that there are at least two minimal
 normal subgroups.
 Let $M_1,\, M_2$ be distinct minimal normal subgroups of $G$.
 Then $M_1 \cap M_2 = \{1\}$ so $M_1,\, M_2$ centralise 
 each other and generate their direct product.
 Moreover since each is transitive on $\Omega$,
 $M_1$ is the full centraliser of $M_2$ in 
 $\Sym(\Omega)$ (and vice-versa) and hence 
 $M_1, M_2$ are the only  minimal normal subgroups---%
 thus $G$ is of bilithic type. 
 Each acts regularly on $\Omega$, and so if 
 $$
  M_0 := H \cap (M_1 \times M_2)= (M_1 \times M_2)_\alpha
 $$ 
 then $M_1 \times M_2 = M_0.M_1 = M_0.M_2$.
 Thus $M_0$ projects isomorphically onto each of $M_1$ and $M_2$, 
 and $M_1 \cong M_2 \cong M_0$. 
 
 Let $L$ be a non-trivial normal subgroup of $H$ contained 
 in $M_0$ and let $L_1,\,L_2$ be the projections of $L$ 
 into $M_1,\, M_2$ respectively. 
 Since $L \Normal M_0$ we have $L_i \Normal M_i$.
 Now $H$ also normalises  $L_i$, and it follows
 that $L_i \Normal G$.
 Therefore by minimality of $M_i$ we have $L_i = M_i$ 
 and so $L = M_0$.
 Thus $M_0$ is a minimal normal subgroup of $H$.

 This completes the proof Theorem~$\ref{thm:maxfin}$.
 
 \bigskip
 Now suppose that $H$ is finite.
 There exist non-trivial normal subgroups of $M$
 whose distinct $H$-conjugates are pairwise 
 disjoint---rather trivially, for example, $M$
 itself satisfies this condition.
 Choose $T \Normal M$ (and $T \neq \{1\}$) such that
 the distinct $H$-conjugates  $T_1,\, T_2,\, \ldots,\, T_q$ 
 of $T$ are pairwise disjoint and furthermore
 $q$ is as large as possible subject to this condition.
 By construction, if $i \neq j$ then $T_i \cap T_j = \{1\}$ 
 so $T_i,\; T_j$, being normal subgroups of $M$,
 commute elementwise.
 In particular, $T_1$ centralises $T_2T_3\,\cdots\,T_q$.
 Let $Z_1 := T_1 \cap (T_2T_3\,\cdots\,T_q)$.
 Then $Z_1 \leq Z(T_1)$, so $Z_1$ is abelian.
 Now $Z_1 \Normal M$ so $Z_1$ has only finitely many
 conjugates $Z_1,\, Z_2,\, \ldots,\, Z_s$ in $G$
 and their product $Z_1Z_2\, \cdots\, Z_s$ is 
 normal in $G$.
 Therefore either $Z_1 = \{1\}$  
 or $Z_1Z_2\, \cdots\, Z_s = M$.

 Suppose (seeking a contradiction) that  $Z_1 \neq \{1\}$.
 Now, being a product of abelian normal subgroups
 (of itself), by Fitting's Theorem $M$ is nilpotent and its centre 
 $Z$ is a non-trivial abelian normal subgroup of~$G$.
 By minimality $M = Z$, that is, $M$ is abelian.
 For $a \in M\setminus\{1\}$ the group 
 $\langle h^{-1}ah \mid h \in H\rangle$
 is normalised by both $H$ and $M$,
 and is therefore a non-trivial normal subgroup of $G$,
 hence is equal to~$M$.
 Thus $M$ is finitely-generated and abelian. 
 But that is impossible: since there are no non-trivial
 finite normal subgroups in $G$, $M$ must be free 
 abelian of finite rank and so $\{a^2 \mid a \in M \}$ 
 is a normal subgroup of $G$ properly contained in~$M$. 
 This contradiction shows that $Z_1 = \{1\}$, that is, that 
 $T_1 \cap (T_2T_3\,\cdots\,T_q) = \{1\}$.
 
 Similarly, of course $T_i$ intersects the product of
 the groups $T_j$ for $j \neq i$ trivially.
 Therefore $M = T_1 \times T_2 \times \cdots \times T_q$
 and $H$ acts by conjugation to permute the factors 
 $T_i$ transitively.
 Let $U_1$ be a non-trivial normal subgroup of $T_1$ 
 (so $U_1\Normal M$).
 If $h \in H$ then $U_1^h \Normal T_1^h$, so $U_1^h \Normal M$ and
 \hbox{$\langle U_1^h \mid h \in H\rangle \Normal G$}.
 Since $U_1$ has a conjugate inside each of 
 $T_1,\, T_2,\, \ldots,\, T_q$, there are at least $q$ 
 conjugates of $U_1$ in $M$ which are pairwise disjoint.
 By the maximality of $q$ therefore, $U_1$ has exactly one
 conjugate $U_i$ in $T_i$ for each $i \in [1\,..\,q]$
 and, arguing as above, $M = \langle U_1^h \mid h \in H\rangle 
 = U_1 \times U_2 \times \cdots \times U_q$.
 It follows that $U_1 = T_1$, and so $T_1$ is simple. 
 This completes our proof of Theorem~\ref{thm:finitestabs}.\\

 A similar argument may be used to show the following. 

 \begin{obs}\label{Obs:minSN}
 If $G$ is normally infinite of bilithic type
 and $H$ satisfies \minSN\ then $M_0$, $M_1$ and $M_2$
 {\rm(}as in the statement of Theorem~{\rm\ref{thm:maxfin})}
 are direct products of finitely many isomorphic 
 infinite simple groups.
 \end{obs}
 
 For, normal subgroups of $M_0$ are subnormal in $H$,
 and therefore if $H$ satisfies \minSN\ then $M_0$
 satisfies \minN.
 Let $T$ be a minimal normal subgroup of $M_0$.
 For a finite subset $\Phi$ of $H$ define
 $P_\Phi := \langle T^h \mid h \in \Phi\rangle$ and 
 $C_\Phi := C_{M_0}(P_\Phi)$.
 Then $P_\Phi \Normal M_0$ and so $C_\Phi \Normal M_0$.
 Since $M_0$ satisfies \minN, there exists a finite
 subset $\Psi$ of $H$ such that $C_\Psi$ is minimal 
 in the set $\{C_\Phi \mid \Phi \subseteq_{\mathrm{fin}} H\}$.
 Then $C_\Psi = C_{\Psi \cup \{h\}}$ for any $h \in H$,
 and so $C_\Psi = C_{M_0}(P)$ where 
 $P := \langle T^h \mid h \in H\rangle$.
 Clearly, $P \Normal H$ and $P \leq M_0$ and so 
 $P = M_0$ since $M_0$ is a minimal normal subgroup of $H$.
 Then $C_\Psi = C_{M_0}(P) = C_{M_0}(M_0)$.
 If $C_\Psi = M_0$ then $M_0$ is abelian and 
 since it satisfies \minN, it would have to be finite, 
 which is not the case.
 Therefore $C_\Psi = \{1\}$ (being the centraliser of $P$ it is 
 normal in $H$).
 Now if $h \in H$ and $T^h \not \leq P_\Psi$ then
 $T^h \cap P_\Psi = \{1\}$ (since $T^h$ is a minimal normal
 subgroup of $M_0$) and so $T^h$ would centralise
 $P_\Psi$ which is not the case.
 Therefore $T^h \leq P_\Psi$ for all $h \in H$, that is,
 $P_\Psi = P$.
 Thus $M_0$ is a product of finitely many conjugates of 
 $T$, and now the proof can be completed as in the 
 case where $H$ is finite.

 \section{Proof of Theorem~\ref{ThmAffine}}\label{s:3} 
 
 Now suppose that there are arbitrarily large finite 
 normal subgroups in~$G$. Let $\calK$ be the set of all 
 finite normal subgroups of $G$ and let 
 $K := \langle N\mid N\in\calK\rangle$.
 Then $K$ is an infinite normal
 subgroup of $G$, hence transitive on~$\Omega$.
 Also let $\calC := \{ C_G(N)\mid N \in \calK\}$. 
 Note that for each $N \in \calK$, $C_G(N)$ is the kernel 
 of the map from $G$ to $\Aut(N)$ induced by conjugation
 and so is normal in $G$.
 Since by Observation~\ref{Obs:GminN}, $G$ satisfies \minN, 
 $\calC$ contains a minimal member $C$. 
 Then $C = C_G(L)$, for some $L \in \calK$.
 Since $L$ is finite, and $G/C \leq\Aut L$,
 $|G:C|$ is finite, so $C$ is infinite, hence transitive. 
 For any $N \in \calK$, since $L$ and $N$ are finite 
 normal subgroups of $G$, so also is $NL$ and hence 
 $NL \in \calK$. 
 Thus $C_G(NL)\in\calC$ so $C \leq C_G(NL)$ by minimality. 
 On the other hand $C_G(NL)$ centralises $L$ and 
 hence $C_G(NL) \leq C$: thus $C_G(NL) = C$. 
 It follows that $C$ centralises  each $N \in \calK$, 
 and hence $C$ centralises $K$. 
 Since both $C$ and $K$ are transitive, we must 
 have $C = C_G(K)$, $K = C_G(C)$, and both $C$ and $K$ 
 are regular on~$\Omega$.

 Let $M := C \cap K$, the centre $Z(K)$. 
 Since $C$ has finite index in $G$, $M$ has finite index
 in $K$, and hence is infinite and therefore transitive. 
 Being a transitive subgroup of the regular 
 groups $C$ and $K$, the group $M$ acts regularly, 
 and it follows that $M = K = C$. 
 Thus $C_G(M) = C_G(K) = C = M$ 
 so $M$ is the unique minimal infinite normal 
 subgroup of $G$. 
 Also, $G = MH$ with $M \cap H=1$, since $M$ is regular,  
 and since $M = C$, it has finite index in~$G$.
 Therefore $H$ is finite and acts faithfully 
 on $M$ by conjugation.

 We next determine the structure of $M$. 
 Since $M$ is abelian we now use additive notation. 
 By the Primary Decomposition Theorem  
 $M = \bigoplus M_p\,$, where the 
 sum is over all prime numbers $p$ and
 $M_p$ is the $p$-primary component of~$M$
 (recall that $M$, being a union of finite normal 
 subgroups is periodic).
 Let $p$ be a prime number such that $M_p \neq \{0\}$.
 If $M_p$ were finite then $\bigoplus_{q \neq p} M_q$
 would be infinite, hence transitive (since it is normal in $G$),
 so we would have $M \leq \bigoplus_{q \neq p} M_q$, which is 
 not the case since $M_p \neq \{0\}$.
 Therefore $M_p$ is infinite.
 As it is a normal subgroup of $G$ it is transitive,
 so $M_p = M$.
 That is, $M$ is a $p$-group.
  
 For positive integers $n$, 
 define $M[n] := \{x \in M \Mid nx = 0\}$.
 We  show next that $M[p]$ is finite. 
 Suppose, seeking a contradiction, that $M[p]$
 is infinite.
 Then it is an infinite normal subgroup of $G$ and so 
 $M[p] = M$, that is, $M$ is elementary abelian of 
 exponent $p$ and infinite rank.
 We consider
 $M$
 as an $A$-module, where $A$ is the 
 finite-dimensional (indeed, finite) algebra $\Field_pH$.
 By Corollary~\ref{Cor:infdim} below, $M$ contains infinite 
 proper submodules (AC is needed here).
 These are infinite normal subgroups of $G$
 that are not transitive, contradicting our assumption.
 This proves that $M[p]$ is finite.

 It follows that $M$ has finite rank $r$ (equal to 
 the rank of the elementary abelian group $M[p]$). 
 Consider the subgroup $pM$, that is $\{px \mid x \in M\}$.
 The map $x \mapsto px$ is an endomorphism $M \to M$
 and its kernel is $M[p]$, which is finite.
 Therefore $pM$ is an infinite subgroup of $M$, 
 obviously characteristic, hence normal in $G$.
 It follows that $pM = M$, that is, that $M$
 is $p$-divisible.
 Being a $p$-group, it is $q$-divisible for all
 prime numbers $q \neq p$, and therefore 
 it is divisible.

 The proof of Theorem~\ref{ThmAffine} is now completed by 
 Theorem~\ref{Thm:H p-adic irreducible} below.

 \section{Some relevant representation theory}\label{s:HM}  

 We begin with a lemma, and a corollary that is needed 
 in the proof 
 of Theorem~\ref{ThmAffine}.
 Recall that the socle $\Soc M$ of a module $M$
 is defined to be the submodule generated by all the simple
 submodules of $M$.
  
 \begin{lem}\label{l:peter} 
 Let $A$ be a finite-dimensional algebra over a 
 field $F$ and let $M$ be an $A$-module. 
 If $\dim(\Soc M)$ is finite then also $\dim M$ is finite. 
 \end{lem}

 {\it Proof}.
 Let $J := \mathrm{Rad}\,A$, the Jacobson radical
 defined as the intersection of all the maximal right 
 ideals of $A$. 
 As is well known, $J$ is nilpotent and annihilates
 any semisimple (right) $A$-module.
 For an $A$-module $M$ define the ascending Loewy series 
 by 
 $$
  L_0 := \{0\},\quad
  L_{i+1}/L_i = \Soc(M/L_i)\kern4pt\mbox{for $i \geq 0$}.
 $$ 
 Since $L_iJ \leq L_{i-1}$ for $i \geq 1$, it follows 
 easily that
 $L_i = \{x \in M \Mid xJ^i = \{0\}\,\}$ for $i \leq k$,
 where $J^k =\{0\}$.
 In particular, $L_k = M$.
 The Loewy length of $M$ is defined to be the smallest $m$
 such that $L_m = M$.

 The assertion of the lemma is trivially true if 
 $M$ is semisimple (Loewy length $1$), so suppose 
 as inductive hypothesis that $m > 1$ and the assertion 
 is known to be true for modules of Loewy 
 length $\leq m - 1$.
 Suppose that the Loewy length of $M$ is $m$
 and $\dim (\Soc M) = n$.
 Let $u_1, \ldots, u_r$ be generators of $J$.
 Consider the map $\mu_i: L_2 \to M$, $x \mapsto xu_i$.
 Since $u_i \in J$ and $L_2/L_1$ is semisimple, 
 $\mathrm{Image}(\mu_i) \leq L_1$.
 Therefore $\mathrm{codim}(\ker \mu_i) \leq n$.
 Now 
 $$
  \Soc M = \{x \in M \mid xJ = 0\}
  =  \ker \mu_1 \cap \ker \mu_2 \cap \cdots \cap \ker \mu_r 
 $$ 
 and therefore $\mathrm{codim}_{L_2}(L_1) \leq r\K n$.
 Thus $\dim(\Soc(M/L_1))$ is finite.
 By
 the
 inductive hypothesis, $\dim(M/L_1)$ is finite,
 and therefore $\dim M$ is finite.

 \begin{rem}
 It is clear that one can derive a bound on 
 $\dim M$ in terms of $\dim(\Soc M)$ and $\dim A$
 from the above argument.
 That bound is unrealistically large, however.
 Using only slightly more sophisticated machinery
 (see \cite[\S\S56, 57, 60]{CR}) 
 we can see that $\dim M \leq \dim A \times \dim(\Soc M)$.
 For, if $\Soc M \cong S_1 \oplus \cdots \oplus S_r$, 
 where the summands $S_i$ are simple, 
 then there is an embedding
 $M \leq U_1 \oplus \cdots \oplus U_r$ where $U_i$ 
 is the injective hull of $S_i$.
 Now the $F$-dual $U^*$ of an injective $A$-module $U$ 
 is a projective module over the opposite algebra
 $A^{\mathrm{op}}$.
 Since $S_i$ is simple, $U_i$ is indecomposable, and so
 $U_i^*$ is also indecomposable and therefore 
 isomorphic to  a summand of the free  
 $A^{\mathrm{op}}$-module of rank~$1$.
 Thus 
 $\dim U_i = \dim U_i^* \leq \dim M^{\mathrm{op}} = \dim M$.
 Therefore $\dim M \leq r.\dim M \leq \dim(\Soc M).\dim M$.    
 \end{rem}
 
 \begin{cor}\label{Cor:infdim}
 If $A$ is a finite-dimensional algebra over a field $F$,
 and $M$ is an infinite-dimensional $A$-module 
 then $M$ has $2^{\aleph_0}$ distinct 
 infinite-dimensional proper submodules. 
 \end{cor}
 
 {\it Proof.}\quad
 By Lemma~\ref{l:peter}, since $A$ is finite-dimensional and
 $M$ is infinite-dimensional also $\Soc M$ is 
 infinite-dimensional.
 Being a sum of simple\break 
 submodules, $\Soc M$ is actually a 
 direct sum of infinitely many simple\break 
 $A$-submodules
 (AC is essential here).
 Therefore $\Soc M$ contains a direct sum
 $\bigoplus_{i\in \Nats} S_i$ of simple $A$-modules
 (the fact that an infinite set contains a countably infinite 
 subset also requires AC, albeit only a weak version). 
 Thus if $I$ is any infinite proper subset
 of $\Nats$ then $\bigoplus_{i \in I} S_i$ is an 
 infinite-dimensional 
 proper submodule, and different choices of $I$ 
 give different submodules.
 Since there are $2^{\aleph_0}$ different possibilities 
 for $I$ there are  $2^{\aleph_0}$ different proper 
 infinite-dimensional submodules of $\Soc M$, 
 hence of~$M$.

 \bigskip
 Next we turn to the analysis of pairs $(H, M)$ where
 $M$ is the divisible
 abelian $p$-group of rank~$r$ that is the minimal
 transitive normal subgroup of~$G$ when $G$ is of 
 $p$-divisible affine type, and $H$, a stabiliser in $G$,
 is finite and acts faithfully on $M$ by conjugation.
 Since any infinite $H$-invariant subgroup of $M$ 
 is normal in $G$, hence transitive, there are no 
 infinite proper $H$-invariant subgroups of
 $M$.
 By the structure theorem for divisible abelian groups
 (see for example \cite[Theorem~23.1]{Fuchs}---AC is 
 required for this),
 $M \cong C_{p^\infty} \oplus 
 \cdots \oplus C_{p^\infty}$, with $r$ summands,
 where $ C_{p^\infty}$ denotes the Pr\"ufer $p$-group 
 (see Example~\ref{Ex:Pruefer}).
 It is not hard to see that such a direct sum 
 decomposition of $M$ leads to 
 an isomorphism of the endomorphism ring $\End M$ with
 the algebra  $\mathrm{M}(r, \hat{\Ints}_p)$
 of $r \times r$ matrices over the ring of $p$-adic integers.
 We may identify $H$ with a subgroup of $\Aut M$, and
 since $\Aut M \cong \GL(r,\hat{\Ints}_p)$, we have an
 embedding $H \leq \GL(r,\hat{\Rats}_p)$, where
 $\hat{\Rats}_p$ is the field of $p$-adic rational numbers.
 
 \begin{thm}\label{Thm:H p-adic irreducible}
 As subgroup of $\GL(r,\hat{\Rats}_p)$, $H$ is irreducible.
 \end{thm}

 {\sc Note}.\quad 
 Under these circumstances we say that $H$ acts
 {\it $p$-adic irreducibly\/} on $M$, or that $H$ is a 
 {\it $p$-adic irreducible\/} group of automorphisms of $M$. 
 This is the definition that completes the statement
 of Theorem~\ref{ThmAffine}.
 
 \bigskip
 {\it Proof of the theorem}.\quad
 Consider the Pontryagin dual $M^*$ of $M$ defined by 
 $M^* := \Hom_\Ints(M, S^1)$, where 
 $S^1 := \{z \in \Complex \Mid |z| = 1\}$.
 Since $M = r.C_{p^\infty}$, 
 $M^* = r.\Hom(C_{p^\infty},C_{p^\infty}) = r.\hat{\Ints}_p$,
 a free $\hat{\Ints}_p$-module of rank $r$.
 Let $W := \hat{\Rats}_p \otimes_{\hat{\Ints}_p} M^*$,
 an $r$-dimensional vector space over $\hat{\Rats}_p$.
 Then $M^* \leq W$ and for every $w \in W$ there exists
 $k \geq 0$ such that $p^k w \in M^*$, and so we may 
 think of $W$ as $p^{-\infty}M^*$.
 Also,
 $$
  H \leq \Aut M = \GL(r,\hat{\Ints}_p) =
  \Aut M^* \leq \Aut W = \GL(r,\hat{\Rats}_p).
 $$
 Let $U$ be a non-zero $H$-invariant subspace of $W$
 and let $s := \dim U$.
 We aim to prove that $U = W$, that is, $s = r$.
 To this end define $U_0 := M^* \cap U$ and 
 $$
  M_0 := U_0^{\kern2pt\perp} := 
  \{x \in M \mid \mbox{$u(x) = 0$ for all $u \in U_0$}\}.
 $$
 As $\hat{\Ints}_p$-modules, 
 $M^*/U_0 = M^*/(U \cap M^*) \cong (U + M^*)/U \leq W/U$,
 so $M^*/U_0$ is torsion-free.  
 It is also finitely generated.
 Therefore $M^*/U_0$ is free since $\hat{\Ints}_p$ is a 
 principal ideal domain, and so 
 $M^* \cong U_0 \oplus (M^*/U_0)$.
 Thus $U_0$ is a free summand of $M^*$; it is of rank $s$ 
 since if $u_1,\, \ldots,\, u_s$ is a basis for $U$ 
 then there exists $k \in \Nats$ such that 
 $p^k u_i \in M^*$ for each $i$ and these $s$
 elements are $\hat{\Ints}_p$-independent. 
 Clearly it is $H$-invariant.
 It follows easily that $M_0$ is an $H$-invariant 
 summand of $M$ of rank $r - s$.
 There are no infinite proper $H$-invariant subgroups 
 of $M$ and so, since $s \geq 1$, it follows that 
 $s = r$ and $U = W$.  
 Thus $H$ is an irreducible subgroup of 
 $\GL(r,\hat{\Rats}_p)$, as required.

 \bigskip
 We turn now to the pair $(H,V)$, where
 $V = M[p]$ construed as an $\Field_pH$-module 
 of dimension $r$.
 Earlier we had erroneously persuaded ourselves that
 $V$ must be irreducible as 
 $\Field_pH$-module.
 That need not be true, as is shown by the following
 example that we owe to Peter Kropholler and Karin Erdmann.
 
 \begin{EXA}\label{ex:Cfour}
 \rm 
 The group $H$ generated by the matrix 
 $\begin{pmatrix} 
 \phantom{-}0 &1\\ -1&0
 \end{pmatrix}$ 
 over $\hat{\Ints}_2$ is cyclic of order $4$ and 
 irreducible over $\hat{\Rats}_2$, and so the
 split extension of $C_{2^\infty} \oplus C_{2^\infty}$
 by $H$ has an action of $2$-divisible affine type.
 In this case the action of $H$ on $V$ has kernel 
 of order $2$ and $H$ acts reducibly on $V$ as a 
 cyclic group of order $2$.
 \end{EXA} 
 
 That $H$ need not act faithfully on $V$ is shown 
 already by the simpler example of the generalised
 dihedral group $G := D_{2^\infty}$, the split extension 
 of the Pr\"ufer $2$-group by a cyclic group of 
 order $2$ whose generator acts as inversion.
 Our next example shows that the kernel $K$ of 
 the action of $H$ on $V$ can be arbitrarily large. 
 
 \begin{EXA}\label{ex:LargeKernel}
 \rm
 Let $L := C_{2^\infty} \oplus C_{2^\infty}$
 with the action of $C_4$ described in 
 Example~\ref{ex:Cfour}.
 Let $r := 2s$ where $s \geq 2$, let $M := sL$ 
 (the direct sum of $s$ groups, each isomorphic to $L$),
 and let $H := C_4 \wreath \Sym(s)$. 
 The natural imprimitive action of the wreath product 
 $H$ on $M$ is faithful and $2$-adic irreducible. 
 The kernel of the $H$-action on $M[2]$ is $K$, where
 $K := C_2^{\kern2pt s} \leq C_4^{\kern2pt s}$.
 Thus in this example $K$ is an elementary abelian
 $2$-group of order~$2^s$.
 \end{EXA} 

 \begin{thm}\label{Th:faithful} 
 If $p > 2$ then $H$ acts faithfully on $V$.
 If $p = 2$ then $K$, the kernel of the action of $H$ on $V$,
  is an elementary abelian $2$-group.
 \end{thm} 
 
 {\it Proof}.\quad
 Let $a \in K$. 
 Consider $a - 1 \in \End(M)$.
 Since $a - 1$ annihilates $M[p]$ and the annihilator 
 of $M[p]$ in the endomorphism ring $\End(M)$ 
 is $p\K \End(M)$, we may write
 $a = 1 + pX$ for some $X \in \End(M)$. 
 Suppose now that $a \neq 1$ and (without loss of generality)
 that $a$ has prime order $q$. 
 Then $X \neq 0$ and so there is a non-negative integer 
 $v$ such that $X \in p^v\K\End(M)$, $X\notin p^{v+1}\K\End(M)$, 
 that is, $v$ is the $p$-adic valuation $v_p(X)$ of $X$. 

 If $q \neq p$ then $1 = (1 + pX)^q \equiv 1 + qpX
 \pmod{p^{2v+2}\K\End(M)}$, whence $qpX \in p^{2v+2}\K\End(M)$, 
 which is not the case.  
 Hence $q = p$, and it follows that $K$ is a $p$-group.
 Next suppose that $p$ is odd. 
 Then 
 \begin{eqnarray*}
 1 = (1 + pX)^p 
 &=& 
 1 +p^2X + \binom{p}{2} p^2X^2 + \cdots + p^pX^p\\
 &\equiv& 
 1 + p^2X \pmod{p^{2v+3}\K\End(M)}.
 \end{eqnarray*}
 This implies that $p^2X\in p^{2v+3}\End(M)$, 
 which is false since $v_p(p^2X) = v + 2$. 
 Therefore $K = \{1\}$ if $p$ is odd. 

 Suppose now that $p=2$, and that $a \in K$ has order $4$. 
 Then $a = 1 + 2X$ where $X \neq 0$.
 Since $a^2 = (1 + 2X)^2 = 1 + 4X + 4X^2 \neq 1$ it follows that
 if $Y := X + X^2$ then $Y \neq 0$. 
 Let $w := v_2(Y)$. 
 Now $a^2 = 1 + 4Y \neq 1$ and 
 $$
 1 = a^4 = (1 + 4Y)^2 = 1 + 8Y + 16Y^2.
 $$
 Thus $8Y + 16Y^2 = 0$. 
 Since $v_2 (8Y) = w + 3$ while $v_2(16Y^2) = 2w + 4$,
 however, this is impossible. 
 Thus $K$ is of exponent dividing $2$
 and is an elementary abelian $2$-group,
 as in the statement of the theorem.

 \bigskip
 Now begin with a prime number $p$ and a pair $(H,V)$, 
 where $H$ is a finite group and $V$ is an 
 $\Field_pH$-module of dimension $r$.
 If there exists a divisible abelian $p$-group $A$ of rank~$r$ 
 and an embedding $H \leq \Aut A$ such that 
 $A[p] \cong V$ as $\Field_pH$-module then we call 
 $A$ a {\it divisible hull\/} of $V$ and write $A = p^{-\infty}V$.
 In this language the question to be addressed is: 
 \begin{center}
 \parbox{280pt}
 {what conditions on the pair $(H,V)$ ensure the 
 existence of a divisible hull $p^{-\infty}V$?}
 \end{center}
 
 \noindent
 By Theorem~\ref{Th:faithful} it is necessary that $H$ 
 acts faithfully on $V$ if $p$ is odd and that if $p = 2$ 
 then the kernel of the action is an elementary abelian 
 $2$-group.
 This condition is very far from sufficient, however,
 as Example~\ref{Ex:nocando} below~shows.
 
 Let $M$ be an $RH$-module that is $R$-free of rank $r$ and
 which is such that $M/pM \cong_H V$, where $R$ is some 
 integral domain of characteristic $0$ such
 that $R/pR \cong \Field_p$.
 We call $M$ an {\it integral cover\/} of $V$
 (in the literature it is also called an $R$-form, 
 but since we do not wish to specify $R$, we prefer a less
 specific term).
 The following lemma will prove~useful.
 
 \begin{lem}\label{lem:dual}
 If $V$ has an integral cover then also $V^*$, 
 the dual $\Field_pH$-module, has an integral cover.
 \end{lem}
 
 For, if $M$ is an integral cover of $V$ and 
 $M^* := \Hom_R(M,R)$, where $R$ is the relevant integral domain,
 then $M^*$ is also a free $R$-module, and of the same rank $r$.
 The natural map $R \to R/pR = \Field_p$ induces a 
 homomorphism $M^* \to \Hom_R(M, \Field_p)$ with 
 kernel $pM^*$.
 Every member of $\Hom_R(M, \Field_p)$ has $pM$ in its kernel,
 and so there is a natural isomorphism 
 $\Hom_R(M, \Field_p) \cong \Hom_{\Field_p}(M/pM, \Field_p)
 \cong \Hom_{\Field_p}(V, \Field_p) = V^*$.
 That is, reduction modulo $p$ provides an isomorphism
 $M^*/pM^* \to V^*$.
 Therefore $V^*$ has $M^*$ as an integral cover.

 \bigskip
 In general $V$ need not have either a divisible hull
 or an integral cover. 
 The two go together, however:
 
 \begin{thm}\label{Th:integralcover}
 The finite-dimensional $\Field_pH$-module $V$
 has a divisible hull if and only if it has 
 an integral cover.
 \end{thm}

 {\it Proof.}\quad
 Suppose first that $V$ has an integral cover $M$,
 an $RH$-module for some integral domain $R$ of 
 characteristic $0$ with $R/pR \cong \Field_p$.
 Let $F$ be the field of fractions of $R$ and let
 $$
  S := p^{-\infty}R := \{a/p^k \mid a \in R,\; k \in \Nats\}
  \subseteq F\,.
 $$
 Then $S$ is a subring of $F$ and $R \leq S$.
 Define $p^{-\infty}M := S\otimes_R M$.
 Since $M$ is a free $R$-module of rank $r$, $p^{-\infty}M$
 is an $SM$-module that is free of rank $r$ as $S$-module .
 It contains $M$ as an $RH$-submodule, and 
 $p^{-\infty}M/pM$ is an $RH$-module $A$ with the property 
 that $A[p] \cong M/pM \cong V$ as $\Field_pH$-module.
 Thus $V$ has a divisible hull.
 
 Now suppose conversely that $V$ has a divisible hull $A$.
 Consider the dual group $A^* := \Hom(A, C_{p^\infty})$.
 As in the proof of Theorem~\ref{Thm:H p-adic irreducible},
 $A^*$ is an $RH$-module where $R = \hat{\Ints}_p$
 and $A^*$ is $R$-free of rank $r$.
 Each element $\varphi \in A^*$ induces a homomorphism
 $A[p] \to C_{p^\infty}[p]$ and so there is a restriction map
 $\rho: A^* \to \Hom(A[p], C_{p})$ (where $C_p$
 denotes the cyclic group of order $p$). 
 It is not hard to see that
 $
 \ker \rho = \{\varphi: A \to C_{p^\infty} \mid 
 A[p]\leq \ker \varphi\} = pA^*
 $. 
 Therefore $A^*/pA^* \cong \mathrm{Image}(\rho) = 
 \Hom(A[p], C_p) = V^*$.
 Thus $A^*$ is an integral cover of $V^*$.
 Since $V^{**} = V$, applying Lemma~\ref{lem:dual} 
 to $V^*$ we see that $V$ has an integral cover, 
 as required. 
 
 \bigskip
 Finite groups $H$ with $\Field_pH$-modules $V$ that 
 have no integral cover (and therefore no 
 $p$-divisible hull) certainly exist:
 
 \begin{EXA}\label{Ex:nocando}
 \rm 
 If $p \geq 5$, $H := \GL(2,p)$ and $V$ is the natural
 $2$-dimen\-sional module $\Field_p^{\kern2pt2}$ then 
 $V$ has no integral cover.
 For, if $R$ were an integral domain of characteristic $0$
 with field of fractions $F$, 
 and $M$ an $RH$-module that is $R$-free of rank $2$,
 then $F \otimes_R M$ would be an $FH$-module of 
 dimension~$2$ with $H$ acting faithfully.
 But $H$ has a subgroup isomorphic to the metacyclic
 group $\AGL(1,p)$ and it is easy to see 
 that this has no faithful representation of dimension
 $< p - 1$ over any field of characteristic $\neq p$.
 Therefore $H$ has no faithful representation of 
 degree $< p - 1$ over~$F$.
 \end{EXA}

 {\sc Comment~1}.\quad
 Let us say that $V$ has a $\Ints_{p^2}$-hull $p^{-1}V$ if
 there exists a $\Ints_{p^2}H$-module $X$ that is
 free of rank $r$ as $\Ints_{p^2}$-module and such that 
 \hbox{$X[p] \cong V$} as $\Field_pH$-module.
 The map $x \mapsto px$ will then be an endomorphism 
 of $X$ with kernel and image both isomorphic to $V$.
 Define a $\Ints_{p^k}$-hull analogously.
 If $Y$ were a $\Ints_{p^3}$-hull then $pY$ and $Y[p^2]$
 would be `overlapping'  $\Ints_{p^2}$-hulls.
 Intuition suggests that if a $\Ints_{p^2}$-hull
 exists then one should be able to manufacture a 
 $\Ints_{p^3}$-hull from two overlapping copies; 
 then, by some sort of boot-strapping, a
 $\Ints_{p^k}$-hull for every $k \geq 2$.
 It should follow that $V$ has a divisible hull
 $p^{-\infty}V$ if and only if it has a 
 $\Ints_{p^2}$-hull.
 Is this true?
 
 \bigskip
 {\sc Comment~2}.\quad
 Let $H$ be any finite group and $V$ an 
 $\Field_pH$-module.
 It is not hard to see from a combination of 
 Theorems~\ref{Thm:H p-adic irreducible}
 and~\ref{Th:integralcover} that the pair $(H,V)$
 arises from a permutation group of $p$-divisible affine
 type if and only if $V$ has an integral cover
 over some integral domain $R$ of characteristic $0$
 (not necessarily $\hat{\Ints}_p$) which is rationally
 irreducible in the sense that it is irreducible
 as $FH$-module where $F$ is the field of fractions of $R$.  
 Consider the case that $V$ is irreducible.
 From the beginnings of modular representation theory 
 we see that if $V$ lies in a $p$-block of defect $0$ 
 (in the sense that its constituents when $\Field_p$ 
 is extended to a splitting field lie in blocks of 
 defect~$0$) then $V$ has a rationally irreducible 
 integral cover (or equivalently a
 $p$-adic irreducible integral cover), 
 and therefore $(H,V)$ can arise from 
 a group $G$ of $p$-divisible affine type as 
 in Theorem~\ref{ThmAffine}.
 We had hoped that this condition would be necessary 
 as well as sufficient but that is not the case.
 We are grateful to Karin Erdmann for 
 drawing our attention to examples due to Gordon James 
 (see \cite{James} or \cite[Theorem~7.3.23, Example~7.3.26]{JK}) 
 of modules of non-zero defect that have integral covers.  

 \bigskip
 {\bf Acknowledgements}. 
 The first two authors are grateful to the Queen's College,
 Oxford, All Souls College, Oxford, and 
 the Mathematical Institute, Oxford, for
 hospitality that facilitated their part of the 
 collaboration. 
 The third author was supported by an 
 Australian Research Council Discovery Early 
 Career Researcher Award 
 (DECRA), project number DE130101521.
 It is a pleasure to record our gratitude to 
 Karin Erdmann and Peter Kropholler for helpful 
 criticisms of an earlier draft of this article.

\end{document}